\documentclass[12pt]{article}
\usepackage[T2A]{fontenc}
\usepackage[utf8]{inputenc}
\usepackage{graphicx,xcolor}
\usepackage{amssymb,amsmath,amsfonts,amsthm,amscd,latexsym,verbatim,graphics,epsfig,indentfirst,color}
\usepackage{geometry}
\begin{document}

\begin{center}
    {\Large  Right representations of Novikov algebras}

    A. S. Panasenko
\end{center}

\begin{abstract}
We study the theory of right representations of Novikov algebras. We introduce the notion of a modular right ideal of a Novikov algebra and prove that every irreducible right representation of a Novikov algebra is isomorphic to the quotient by a maximal modular right ideal.

We prove that the class of Novikov algebras with no irreducible representations forms a hereditary radical. We introduce the notion of a primitive Novikov algebra and show that a Novikov algebra is primitive if and only if it admits an almost faithful irreducible representation. We also introduce the notion of the quasi-kernel of a representation, defined as the largest ideal contained in its kernel. Finally, we prove that the constructed radical coincides with the intersection of the quasi-kernels of all its irreducible representations.

\end{abstract}

\section{Introduction}

\textbf{Definition.} An algebra is called a \textbf{Novikov algebra} if it satisfies the identities
\begin{gather}
    (x,y,z)=(y,x,z),\label{LAs}\\ 
    (xy)z=(xz)y, \label{RCom}
\end{gather}
where $(x,y,z)=(xy)z-x(yz)$ denotes the associator.

Novikov algebras were introduced by I.~M.~Gel'fand and I.~Ya.~Dorfman \cite{GD1979} in the context of Hamiltonian operators in the formal calculus of variations. Independently, these algebras emerged in the studies of A.~A.~Balinskii and S.~P.~Novikov \cite{BN1985} (see also \cite{N1985}) concerning local translation-invariant Lie algebras.

The study of finite-dimensional Novikov algebras was initiated by E.~I.~Zel'manov \cite{Z1987}, who classified simple and semisimple finite-dimensional Novikov algebras over algebraically closed fields of characteristic~0. He demonstrated that every such algebra is commutative and, consequently, decomposes into a direct sum of fields. For fields of positive characteristic, the classification of simple finite-dimensional Novikov algebras over algebraically closed fields of characteristic $p>2$ was achieved by J.~M.~Osborn \cite{Osborn1992} and X.~Xu \cite{Xu1996}, with further recent contributions by V. N. Zhelyabin and A. S. Zakharov \cite{ZZ2024}. The classification over arbitrary fields of characteristic $p>0$ was recently completed by V. N. Zhelyabin and A. P. Pozhidaev \cite{PZ2026}.

Regarding the representation theory of these algebras, bimodules over simple Novikov algebras of characteristic zero possessing an idempotent were initially investigated by J. M. Osborn \cite{Osborn1995}, and their complete classification was later provided by X. Xu \cite{Xu2001}. Bimodules over finite-dimensional simple Novikov algebras were also studied by E.~I.~Zel'manov and J.~M.~Osborn \cite{OZ1995}. X.~Xu \cite{Xu1996} classified finite-dimensional irreducible bimodules over finite-dimensional simple Novikov algebras defined over algebraically closed fields of characteristic $p>0$.

Despite these extensive studies on bimodules, the theory of one-sided (right and left) modules over Novikov algebras remains largely unexplored. The only notable precursor is the work of A.~M.~Slinko and I.~P.~Shestakov \cite{SS1974}, who introduced the general notion of right representations for an arbitrary variety of algebras and developed the corresponding theory specifically for alternative and right alternative algebras.

In a series of previous papers, the author investigated the radical theory of Novikov algebras. Specifically, the Baer radical was introduced and analyzed in \cite{P2022}, the Andrunakievich radical was characterized in \cite{P2024}, and the radicals of Lie-solvable Novikov algebras were explored in \cite{P2026}. 

The present paper is dedicated to the systematic study of right representations of Novikov algebras. Our main contributions are as follows. First, we prove that the class of Novikov algebras lacking irreducible right representations forms a hereditary radical. Second, we demonstrate that this radical coincides with the intersection of the quasi-kernels (the largest ideals contained in the kernels) of all irreducible representations, as well as with the intersection of all primitive ideals of the algebra. To achieve this, we introduce the concepts of a modular right ideal and a primitive Novikov algebra, establishing that an algebra is primitive if and only if it admits an almost faithful irreducible representation.







\section{Preliminaries}

Throughout the paper, $F$ denotes a field of characteristic different from $2$, and all algebras are assumed to be defined over $F$.

In every Novikov algebra $A$, the following identities hold \cite{SF2002}:
\begin{gather}
	(x,y,z)t = (xt,y,z) = (x,yt,z) \label{AsLin}, \\
    x[y,z]=(x,z,y)-(x,y,z) \label{RAs}
\end{gather}
for all $x,y,z,t\in A$, where $[a,b]=ab-ba$ denotes the commutator of $a$ and $b$.

Furthermore, the following identity holds in any Novikov algebra \cite{SF2002}:
\begin{equation}
[z,y]x=\frac{1}{2}[z,yx]+\frac{1}{2}[zx,y].\label{AssDiff}
\end{equation}
From this equality it follows that for any ideal $I$, the space $[I,I]$ spanned by all commutators of elements of $I$ is an ideal of the algebra.

\medskip We begin by recalling the definition of a right representation introduced by A.~M.~Slinko and I.~P.~Shestakov \cite{SS1974}, restricting their original theory to algebras over the field $F$, since this theory was originally created for algebras over commutative rings. Let $T(V)$ denote the tensor algebra of a vector space $V$:
\[T(V)=F \oplus V^{(1)}\oplus\dots\oplus V^{(n)}\oplus\dots,\]
where $V^{(n)} = V\otimes \dots \otimes V$ is the $n$-fold tensor product of $V$. Every linear map $\psi:V\to W$ uniquely extends to an algebra homomorphism
\[T(\psi):T(V)\to T(W)\] 
determined by the conditions 
\[T(\psi)(1)=1, \qquad T(v_1\otimes\dots\otimes v_n)=\psi(v_1)\otimes \dots\otimes \psi(v_n).\] 
Consequently, if $B$ is an associative unital algebra, then every linear map $\psi:V \to B$ induces an algebra homomorphism
\[T(\psi): T(V)\to B\]
such that
\[T(\psi)(1)=1,\qquad T(\psi)(v_1\otimes\dots\otimes v_n)=\psi(v_1)\dots\psi(v_n).\]

Let $\mathcal{M}$ be a variety of algebras over a field $F$, and let $\mathcal{M}_X$ denote the free algebra in $\mathcal{M}$ on the countable set of generators $X=\{x_1,x_2,\dots\}$. Let $F(x_1,\dots,x_n)\in T(\mathcal{M}_X)$. If $A$ is an algebra in the variety $\mathcal{M}$, then, for $a_1,\dots,a_n\in A$, we denote by $F(a_1,\dots,a_n)$ the element of $T(A)$ obtained from $F(x_1,\dots,x_n)$ by substituting $a_i$ for $x_i$.

\medskip Let $\varphi:\mathcal{M}_X\to A$ be the algebra homomorphism determined by $\varphi(x_i)=a_i$ for $i=1,\dots,n$. Then 
\[F(a_1,\dots,a_n)=T(\varphi)(F(x_1,\dots,x_n)).\]

\textbf{Definition.} Let $\rho:A\to \mathrm{End}(V)$ be a linear map from an algebra $A$ to the endomorphism algebra of a vector space $V$, and let $T(\rho): T(A)\to\mathrm{End}(V)$ be the induced algebra homomorphism. An element $F(x_1,\dots,x_n)\in T(\mathcal{M}_X)$ is called an \textbf{identity of the map} $\rho$ if $T(\rho)(F(a_1,\dots,a_n))=0$ for all $a_1,\dots,a_n\in A$.

\medskip We now define 
\[R:\mathcal{M}_X\to \mathrm{End}(\mathcal{M}_X)\]
to be the linear map given by $R(f)=R_f$, where $R_f$ denotes the operator of right multiplication by $f\in\mathcal{M}_X$. Let $T(R): T(\mathcal{M}_X)\to \mathrm{End}(\mathcal{M}_X)$ be the algebra homomorphism induced by $R$. We denote the kernel of this homomorphism by $I_{\mathcal{M}}$. 

\medskip\textbf{Definition.} The elements of the ideal $I_{\mathcal{M}}$ are called \textbf{R-identities} of the variety $\mathcal{M}$.

\medskip Let $A$ be the algebra in the variety $\mathcal{M}$ and let $M$ be a vector space. 

\medskip \textbf{Definition.} A linear map $\rho:A\to \mathrm{End}(M)$ is called a \textbf{right $\mathcal{M}$-representation} of $A$ if every R-identity of the variety $\mathcal{M}$ is an identity of the map $\rho$. In this case, the space $M$, equipped with the action $m\cdot a=m\rho(a)$, for $m\in M$ and $a\in A$, is called a \textbf{right $A$-module in the variety $\mathcal{M}$}.

\medskip Left representations and left modules are defined similarly. Since the variety of Novikov algebras is not closed under anti-isomorphisms, left and right representations may have significantly different properties. Throughout the remainder of the paper, we consider only right representations and right modules and therefore omit the adjective \textit{right}.

\medskip The kernel of a representation $\rho$ of an algebra $A$ is denoted by $\mathrm{Ker}_{\rho}(A)$.

\medskip \textbf{Definition.} The largest ideal of $A$ contained in $\mathrm{Ker}_{\rho}(A)$ is denoted by $K_{\rho}(A)$ and is called the \textbf{quasi-kernel} of $\rho$. A representation is called \textbf{faithful} if $\mathrm{Ker}_{\rho}(A)=0$, and \textbf{almost faithful} if $K_{\rho}(A)=0$. An $A$-module $M$ is called \textbf{irreducible} if $M$ has no nonzero proper submodules and $MA\neq 0$. The corresponding representation is also called irreducible.

\medskip The concepts of prime and semiprime algebras play a fundamental role in the structural theory of algebras.

\medskip\textbf{Definition.} An algebra $A$ is called \textbf{prime} if, for any ideals $I$ and $J$ of $A$, the equality $IJ=0$ implies that $I=0$ or $J=0$. An algebra $A$ is called \textbf{semiprime} if, for every ideal $I$ of $A$, the equality $I^2=0$ implies that $I=0$.

\section{Irreducible Novikov Representations}

It follows from the identities \eqref{RCom}--\eqref{RAs} that the ideal of R-identities of the variety of Novikov algebras contains the following elements:
\begin{gather*}
	x_1\otimes x_2 - x_2\otimes x_1,\\
	x_1\otimes x_2\otimes x_3 - (x_1x_2)\otimes x_3 - x_3\otimes x_1 \otimes x_2 + x_3\otimes (x_1x_2),\\
	x_1\otimes x_2\otimes x_3 - (x_1x_2)\otimes x_3 - (x_1x_3)\otimes x_2 + (x_1x_3)x_2.
\end{gather*}

Hence, every Novikov $A$-module $M$ satisfies the following identities for all $m\in M$ and $a,b,c\in A$:
\begin{gather}
	(ma)b=(mb)a,\label{ModRCom}\\
	(m,a,b)c=(mc,a,b)=(m,ac,b),\label{ModAsLin}\\
	m[a,b]=(m,b,a)-(m,a,b).\label{ModRAs}
\end{gather}

\medskip In \cite{SS1974}, the notions of left and right associativity of a module were introduced. For right-commutative algebras, however, these notions coincide. 

\textbf{Definition.} A Novikov $A$-module $M$ is called \textbf{associative} if $(m,a,b)=0$ for all $m\in M$ and $a,b\in A$.

Let $M$ be a Novikov $A$-module and let $B$ a subspace of $A$. Then the subspace $MB=\{\sum m_ib_i\mid m_i\in M, b_i\in B\}$ is a submodule of $M$.

\medskip In what follows, we will use the fact that an algebra with zero quasi-kernel is prime.

\medskip\textbf{Lemma 3.1.} {\it If $\rho$ is an irreducible representation of a Novikov algebra $A$, then the algebra $A/K_{\rho}(A)$ is prime.}

\medskip\textbf{Proof.} Let $B,C$ be ideals of the algebra $A$ and let $BC\subseteq K_{\rho}(A)$. Suppose that $B\nsubseteq K_{\rho}(A)$ and $C\nsubseteq K_{\rho}(A)$. Let $M$ be the $A$-module corresponding to the representation $\rho$. Since $B\nsubseteq K_{\rho}(A)$, it follows that $B\nsubseteq \mathrm{Ker}_{\rho}(A)$ and $MB\neq 0$. Since $M$ is irreducible, $MB=M$. Similarly, $MC=M$. Hence $(MB)C=MC=M\neq 0$. Then there exist $b\in B$ and $c\in C$ such that $(Mb)c\neq 0$. But 
\[((Mb)c)a \overset{\eqref{ModRCom}}{=} ((Ma)b)c\subseteq (Mb)c,\]
so $(Mb)c$ is a nonzero submodule of $M$ and $(Mb)c=M$. Therefore
\[M = (Mb)c = (((Mb)c)b)c = ((M,b,c)b)c \overset{\eqref{ModAsLin}}{=} (M,b^2c,c) = 0,\]
which is a contradiction. \qed

\medskip It is well known that irreducible modules over an associative algebra can be characterized as quotients of the regular module by maximal right ideals satisfying an additional condition. The following theorem shows that an analogous construction holds for Novikov algebras.

\medskip\textbf{Lemma 3.2.} {\it Let $A$ be a Novikov algebra and let $M$ be an irreducible $A$-module. Then exactly one of the following holds:
	\begin{itemize}
		\item $M$ is associative and isomorphic to the quotient module $A/I$, where $I$ is a maximal two-sided ideal $I$ of $A$ such that $A/I$ is a field;
		\item $M$ is nonassociative and isomorphic to the quotient module $A/I$, where $I$ is a maximal right ideal of $A$, and there exist elements $e,f\in A$ such that $a-(a,e,f)\in I$ for any $a\in A$. Moreover, the element $e$ is determined by the element $f$, and $f$ may be chosen to be any element satisfying $(A,A,f)\nsubseteq I$.
	\end{itemize}
}

\medskip\textbf{Proof.} First, observe that the set $S=\{u\in M\mid uA = 0\}$ is an $A$-submodule of $M$. Since module $M$ is irreducible, it follows that $S=0$.

\medskip 1) Suppose that $(M,A,A)=0$. Then $M$ is associative. In this case, $mA$ is a submodule of $M$ for every $m\in M$. Since $MA=M$, there exists an element $m\in M$ such that $mA=M$. The map $\psi:A\to M$, determined by the rule $\psi(a)=ma$, is a module epimorphism. Hence, the module $M$ is isomorphic to the quotient module $A/I$, where $I$ is a right maximal ideal of $A$. 

\medskip Since $A$ is a Novikov algebra, the subspace $[A,A]$ is an ideal. Moreover, 
\[A[A,A] \overset{\eqref{RAs}}{\subseteq} (A,A,A)\subseteq I,\] 
hence $[A,A]\subseteq K_{\rho}(A)$, where $\rho$ is the right representation corresponding to the module $A/I$. Consequently, the algebra $A/K_{\rho}(A)$ is associative and commutative. It follows that $I=K_{\rho}(A)$. Thus $I$ is a two-sided ideal and the algebra $A/I$ is a field.

\medskip 2) Suppose that $(M,A,A)\neq 0$. Then there exists an $m\in M$ such that $(m,A,A)\neq 0$. Hence, there exists $f\in A$ such that $(m,A,f)\neq 0$. Since 
\[(m,A,f)a \overset{\eqref{ModAsLin}}{=} (m,Aa,f)\subseteq (m,A,f)\] 
for every $a\in A$, it follows that $(m,A,f)$ is a nonzero $A$-submodule of $M$. By the irreducibility of $M$, it follows that $(m,A,f)=M$.

\medskip Consider the map $\psi:A\to M$ determined by the rule $\psi(c)=(m,c,f)$. Observe that
\[\psi(cd)=(m,cd,f)\overset{\eqref{ModAsLin}}{=}(m,c,f)d=\psi(c)d,\]
and therefore $\psi$ is an $A$-module homomorphism. The kernel of this homomorphism $I=\mathrm{Ker}(\psi)$ is a submodule of the regular module $A$. Hence, $I$ is a right ideal of the algebra $A$. Since $\psi(A)=(m,A,f)=M$, by the First Isomorphism Theorem for modules, $M\simeq A/I$. As in the classical case, it is straightforward to verify that $I$ is a maximal right ideal.

\medskip Since $(m,A,f)=M$, there exists an element $e\in A$ such that $m=(m,e,f)$. But then for any $a\in A$ we have
\[(m,a,f)=((m,e,f),a,f) = ((m,e,f)a)f-(m,e,f)(af) \overset{\eqref{ModAsLin}}{=} (m,(e,a,f),f),\]
whence $(m,a-(a,e,f),f) = 0$ by \eqref{LAs} and $a-(a,e,f)\in\mathrm{Ker}(\psi)=I$. \qed

\medskip It follows immediately that every module satisfying the first alternative of Lemma 3.2 is irreducible. We now consider the second alternative.

\medskip\textbf{Lemma 3.3.} {\it Let $A$ be a Novikov algebra and let $I$ be a maximal right ideal of $A$. If there exist elements $e,f\in A$ such that $a-(a,e,f)\in I$ for every $a\in A$, then the quotient module $A/I$ is an irreducible Novikov $A$-module.}

\medskip\textbf{Proof.} Let $I$ be a maximal right ideal of a Novikov algebra $A$ such that there exist elements $e,f\in A$ satisfying $a-(a,e,f)\in I$ for every $a\in A$. We now prove that the module $A/I$ is irreducible. It is clear that $A/I$ has no proper nonzero submodules, since submodules of $A/I$ correspond to right ideals of $A$ containing $I$, and $I$ is maximal. It remains to prove that $(A/I)A\neq 0$. Assume that $(A/I)A=0$. Fix an element $a\in A$. Then $0=((a+I),e,f) = a+I$, so $a\in I$, which contradicts the assumption $I\neq A$. \qed

\medskip In the variety of alternative algebras (defined by the identities $(x,x,y)=(x,y,y)=0$), all irreducible modules are associative. Furthermore, the following result holds for both associative and alternative algebras: the kernel of an irreducible module coincides with its quasi-kernel. In particular, the intersection of the kernels of irreducible modules is an ideal. As the following example demonstrates, none of these statements hold for Novikov algebras.

\medskip\textbf{Example 3.4.} a) $\mathrm{char} F = p > 2$ (\cite{Xu1996}). Let $A=L(y_{-1},y_0,\dots,y_{p^{n}-2})$ be a finite-dimensional simple Novikov algebra over a field $F$ of characteristic $p>2$, with the following multiplication table:
\[y_iy_j=\begin{pmatrix}
	i+j+1\\
	j
\end{pmatrix}y_{i+j}+\delta_{i,-1}\delta_{j,-1}ay_{p^n-2}+\delta_{i,-1}\delta_{j,0}by_{p^n-2}.\]
Consider the subspace $I=L(y_0,y_1,\dots,y_{p^n-2})$. It is clear that $I$ is a maximal right ideal of $A$. Then $(y_{-1},y_{-1},y_1) = ay_{p^n-2}y_{1}-y_{-1}y_0 = -y_{-1}-by_{p^n-2}$. Hence, $y_{-1}-(y_{-1},y_{-1},-y_1) = by_{p^n-2}\in I$. Then $x-(x,y_{-1},-y_1)\in I$ for every $x\in A$, and therefore $A/I$ is an irreducible non-associative Novikov $A$-module by lemma 3.3. If $\rho$ is the corresponding representation, then $K_{\rho}(A)=0$ with $\mathrm{Ker}(\rho) = L(y_1,\dots,y_{p^n-2},y_{-1})\neq 0$. Moreover, $I$ is the unique maximal right ideal.

\medskip b) $\mathrm{char} F = 0$ (\cite{F1989}). Let $A = L(\{x_i\}_{i\ge -1})$ be a simple Novikov algebra over a field $F$ of characteristic $0$ with the following multiplication table:
\[x_ix_j=(j+1)x_{i+j}.\]
Consider the subspace $I=L(\{x_i\}_{i\ge 0})$. It is clear that $I$ is a maximal right ideal of $A$. Then $(x_{-1},x_{-1},x_1)=-2x_{-1}$. Hence $x_{-1}-(x_{-1},x_{-1},\frac{1}{2} x_1) = 0 \in I$. Therefore, $A/I$ is an irreducible non-associative Novikov $A$-module by lemma 3.3. If $\rho$ is the corresponding representation, then $K_{\rho}(A)=0$ with $\mathrm{Ker}(\rho) = L(\{x_i\}_{i\neq 0})\neq 0$. Moreover, $I$ is the unique maximal right ideal.

\medskip Thus, we obtain two observations. First, even a simple non-commutative Novikov algebra admits non-associative modules (unlike, for example, the case of alternative algebras). Second, the intersection of kernels of irreducible representations need not be an ideal, even for a finite-dimensional Novikov algebra (in contrast to the classical cases).

\medskip Let us adapt the concept of an irreducible Novikov module to the case of an associative right-commutative algebra.

\medskip\textbf{Lemma 3.5.} {\it Let $A$ be an associative Novikov algebra. Then every irreducible Novikov $A$-module $M$ is an associative $A$-module.}

\medskip\textbf{Proof.} If $M$ is not an associative $A$-module, then by Lemma 3.2 there exist a maximal right ideal $I$ of $A$ and elements $e,f\in A$ such that $a-(a,e,f)\in I$ for every $a\in A$. Since $(a,e,f)=0$, it follows that $a\in I$ for every $a\in A$. This contradicts the assumption $I\neq A$. \qed

\medskip Recall that a right ideal of an associative algebra $A$ is called \textbf{modular} if there exists $e\in A$ such that $a - ea\in I$ for every $a\in A$. If $A$ is a Novikov algebra, then it satisfies the identity $A[A,A]=0$ by \eqref{RAs}, hence $e[A,A]=0$, whence $[A,A]\subseteq I$. Thus, in an associative Novikov algebra, every modular right ideal is a two-sided ideal whose quotient is a commutative algebra. In connection with this observation and Lemma 3.2, we give the following definition.

\medskip\textbf{Definition.} A right ideal $I\neq A$ of a Novikov algebra $A$ is called \textbf{modular} if one of the following conditions holds: 
\begin{itemize}
\item there exist $e,f\in A$ such that $a-(a,e,f)\in I$ for every $a\in A$;
\item $I$ is a two-sided ideal such that $A/I$ is a commutative algebra with a unit.
\end{itemize}

Thus, the notion of a modular ideal is well defined for both Novikov and associative algebras, and the two definitions agree on their common subclass.

\medskip To avoid further confusion, we need to clarify the ambiguity of the term \textit{maximal modular right ideal}.

\medskip\textbf{Lemma 3.6.} {\it Every modular right ideal $I\neq A$ of a Novikov algebra $A$ is contained in a maximal right ideal of $A$, which is also modular.}

\medskip\textbf{Proof.} Consider the set $X$ of right ideals of $A$ that contain $I$ but do not contain $e$. By Zorn's lemma, $X$ has a maximal element $J$. Let $K$ be a right ideal of $A$ such that $J\subsetneq K$. Then $e\in K$ and, since $(a,e,f)\overset{\eqref{LAs}}{=}(e,a,f)\in K$, it follows that $a - (a,e,f) = a\in K$, i.e. $K=A$. Hence $J$ is a maximal right ideal of $A$. \qed

\medskip Let us notice that the associative definition of modular right ideal does not coincide with Novikov definition.

\medskip\textbf{Example 3.7.} Let $A=(F[x],\circ)$, where $f\circ g = f\frac{dg}{dx}$. Then $I=xF[x]$ is a maximal right ideal of Novikov algebra $A$ and $I$ is not a left ideal of $A$ since $1\circ x = 1\notin I$. A module $A/I$ is irreducible by the same equality. But there is not an $e\in A$ such that $a-e\circ a\in I$ for any $a\in A$. Indeed, $1-e\circ 1 = 1$ for any $e\in A$. But $1-(1,-\frac{1}{2},x^2)=0\in I$, therefore $a-(a,-\frac{1}{2},x^2)\in I$ for any $a\in A$.

\medskip So, we can see a lot of differences of Novikov case in contrast to associative case. Kernel and quasikernel of an irreducible representation does not have to coincide. The notion of modular ideal is different. And there are non-associative irreducible modules even for simple Novikov algebras. 

\section{The Jacobson Radical of Novikov Algebras}

There are several approaches to constructing the Jacobson radical in associative algebras. It is currently unclear whether all approaches will lead to the same radical in Novikov algebras. Some approaches differ in different classes of algebras: for example, the concept of quasiregularity in Jordan algebras differs slightly from the classical one. However, the language of irreducible modules is the same in any variety.

\medskip Denote by $\mathcal{J}$ the class of Novikov algebras that do not have irreducible modules. 

\medskip\textbf{Lemma 4.1.} \textit{The class $\mathcal{J}$ is closed under extensions and homomorphic images, i.e. the following statements hold:
\begin{itemize}
    \item if $A$ is an algebra in $\mathcal{J}$ then $A/I$ is in $\mathcal{J}$ for any ideal $I$ of $A$;
    \item if algebras $I$ and $A/I$ lie in $\mathcal{J}$ for some algebra $A$ and its ideal $I$ then $A$ lies in $\mathcal{J}$ as well.
\end{itemize}} 

\medskip\textbf{Proof.} Let $A_0$ be a Novikov algebra in $\mathcal{J}$ and let $I_0$ be a nonzero ideal of $A_0$. Suppose that $A_0/I_0$ has an irreducible module $M_0$. This means that $A_0/I_0$ has a maximal modular right ideal $\overline{J_0}$. It is easy to see that the full preimage $J_0$ of this right ideal in $A_0$ is a maximal modular right ideal of $A_0$. Thus, the class $\mathcal{J}$ is closed under homomorphic images.

\medskip Now let $J$ be an ideal of $A$ such that $A/J$ and $J$ belong to $\mathcal{J}$. Assume that $A\notin\mathcal{J}$. By Lemma 3.2, two cases are possible.

\medskip 1) The algebra $A$ contains a maximal right ideal $I$ such that $a-(a,e,f)\in I$ for every $a\in A$ and the module $A/I$ is not associative. Then $(A/I,A,A)\neq 0$. Note that $(I+J)/J$ is a right modular ideal of $A/J$. By Lemma 3.6 and Lemma 3.3, it means that $I+J=A$.

\medskip Assume that $(A/I,A,J)\neq 0$. It means that there exists $m\in A/I$ such that $(m,A,J)\neq 0$, whence there exists $y\in J$ such that $(m,A,y)\neq 0$. Since $(m,A,y)$ is a nonzero $A$-submodule in $M$, it follows that $(m,A,y)=M$. Consider the map $\psi:A\to M$ with the condition $\psi(c)=(m,c,y)$. Note that $\mathrm{Ker}(\psi)=I$. The existence of an element $x\in A$ such that $m=(m,x,y)$ implies (as in the proof of Lemma 3.2) that $(m,a,y)=((m,(a,x,y),y),$ whence $a-(a,x,y)\in I$.
\medskip We have $x=x_i+x_j$, where $x_i\in I$ and $x_j\in J$. Let $a\in J$. Then $(a,x_i,y)\in I\cap J$, so $a-(a,x_j,y)\in I\cap J$. Thus, $I\cap J$ is a proper modular right ideal of $J$, a contradiction by Lemma 3.6 and Lemma 3.3.

\medskip Thus, $(A/I,A,J)=0$, that is, $(A,A,J)\subseteq I$. It means that we can assume that $f\in I$ and $e\in J$. In particular, $(J,J,J)\subseteq I\cap J$. It means that the $J$-module $J/(I\cap J)$ is an associative $J$-module. If $J^2\subseteq I$, then $(A,A,A)\subseteq I$, a contradiction. Thus, $J/(I\cap J)$ is a nontrivial associative $J$-module. By the homomorphism theorem, $J/(I\cap J)\simeq (I+J)/I = A/I$, so $A/I$ is a nontrivial associative $J$-module, and $J/I\cap J$ is an irreducible $A$-module.

\medskip Let $K$ be the maximal ideal of the algebra $A$ contained in $I\cap J$. Since the $A$-module $J/I\cap J$ is irreducible, then either $[J,J]\subseteq I$ or $[J,J]+I\cap J = J$. We have $J[J,J]\subseteq I\cap J$. Since $J^2\nsubseteq I$, then $[J,J]\subseteq I\cap J$, that is, $[J,J]\subseteq K$. Thus, the algebra $J/K$ is associative and commutative, so $(I\cap J)/K$ is a two-sided ideal. This means that $I\cap J = K$, that is, $I\cap J$ is a two-sided ideal in $A$. Note that there are no non-zero ideals between $I\cap J$ and $J$ since the module $A/I$ is irreducible. Then $J/(I\cap J)$ is a minimal ideal of $A/(I\cap J)$. Since $J^2$ is not embedded in $I$, $J/(I\cap J)$ is a simple algebra \cite{SZ2020}. Its commutativity was proved above, so $J/(I\cap J)$ is a field. It is a contradiction.

\medskip 2) $A$ has a maximal ideal $I$ such that $A/I$ is a field. If $J\subseteq I$, then $A/I$ is an irreducible $A/J$-module, it is a contradiction. This means that $I+J=A$. Then the algebra $A/I=(I+J)/I$ is isomorphic to the algebra $J/I\cap J$. Thus, as in the previous case, $J/I\cap J$ is an irreducible $J$-module, it is a contradiction. \qed

\medskip We need to give a definition of radical class. Let us fix some class of algebras that is closed under taking homomorphic images and ideals (for the purposes of this paper, we can assume that $\mathcal{M}$ is the class of Novikov algebras). Let $\mathcal{R}$ be a subclass of $\mathcal{M}$. Algebras of $\mathcal{R}$ are called $\mathcal{R}$\textbf{-algebras}. An ideal $I$ of an algebra $A$ is called $\mathcal{R}$\textbf{-ideal} if $I$ belongs (as an algebra) to the class $\mathcal{R}$.

\medskip\textbf{Definition}. A subclass $\mathcal{R}$ of $\mathcal{M}$ is called \textbf{radical} (in the Kurosh-Amitsur sense) if the following conditions are satisfied:

(A) A homomorphic image of an $\mathcal{R}$-algebra in $\mathcal{R}$ is also an $\mathcal{R}$-algebra.

(B) Every algebra $A$ in $\mathcal{M}$ has an $\mathcal{R}$-ideal $\mathcal{R}(A)$ containing all $\mathcal{R}$-ideals of $A$.

(C) If $\mathcal{R}(A)$ is the ideal from condition (B), then there are no nonzero $\mathcal{R}$-ideals in the algebra $A/\mathcal{R}(A)$.

\medskip An algebra $A$ is called $\mathcal{R}$-\textbf{semisimple} if $\mathcal{R}(A)=0$. An algebra $A$ is called $\mathcal{R}$\textbf{-radical} if $\mathcal{R}(A)=A$.

\medskip It is well known that the class $\mathcal{R}$ is radical if and only if conditions (A), (D), and (E) are satisfied, where

(D) The class $\mathcal{R}$ is closed under extensions.

(E) If $\{I_i\}$ is a chain of $\mathcal{R}$-ideals in $A$ whose union coincides with $A$, then $A$ lies in the class~$\mathcal{R}$.

\medskip\textbf{Theorem 4.2.} {\it The class $\mathcal{J}$ is radical in the class of Novikov algebras.}

\medskip\textbf{Proof.} Let us prove condition (E). Let $J_i\in\mathcal{J}$ be a chain of ideals in $A$, where $\bigcup\limits_{i\in\Gamma} J_i = A$. Suppose that $A\notin\mathcal{J}$. By Lemma 3.2, two cases are possible.

1) There exist a maximal right ideal $I$ of $A$ and elements $e,f\in A$ such that $a - (a,e,f)\in I$ for any $a\in A$. Then there exists an $m\in\Gamma$ such that $e,f\in J_m$. Then $a - (a,e,f)\in I\cap J_m$ for any $a\in J_m$, hence $I\cap J_m$ is a modular ideal in $J_m$. By Lemma 3.3 and Lemma 3.6, this means that $I\cap J_m=J_m$, that is, $J_m\subseteq I$. Then $A/I$ is an irreducible $A/J_m$-module, it is a contradiction.

2) There exists a maximal ideal $I$ in $A$, where $A/I$ is a field. Let $e\in A$ be an arbitrary preimage of the identity in $A/I$. Then $e\in J_m$ for some $m$. Let $a,b\in A$. Then
\[[a,b]+I\cap J_m = [ea,b] + I\cap J_m \overset{\eqref{AssDiff}}{=} -[e,ba] + 2[e,b]a + I\cap J_m = 0 + I\cap J_m,\]
that is, $[A,A]\subseteq I\cap J_m$ and the algebra $A/I\cap J_m$ is commutative. Then the algebra $J_m/(I\cap J_m)$ is also commutative and isomorphic to $(I+J_m)/I = A/I$. This means that $J_m/(I\cap J_m)$ is an irreducible $J_m$-module, it is a contradiction.

Conditions (A) and (D) were proved in Lemma 4.1. \qed

\medskip Theorem 4.2 means that in any Novikov algebra $A$, there exists a largest $\mathcal{J}$-ideal $\mathcal{J}(A)$, and the algebra $A/\mathcal{J}(A)$ does not contain nonzero $\mathcal{J}$-ideals. 

\medskip \textbf{Definition.} The radical $\mathcal{J}$ will be called (in accordance with the classical case) the \textbf{Jacobson radical}.

\section{Primitive Novikov Algebras}

There are many different definitions of a primitive algebra, including two key ones: in terms of modular ideals and in terms of irreducible representations. However, the definition in terms of irreducible representations has only been given in varieties where every irreducible almost faithful module is faithful. This is not the case for Novikov algebras, as stated in Example 3.4. Thus, it seems most appropriate to define a primitive algebra in terms of modular right ideals.

\medskip\textbf{Definition.} A Novikov algebra $A$ is called \textbf{primitive} if it has a maximal modular right ideal $I$ that does not contain nonzero two-sided ideals of $A$. An ideal $I$ of $A$ is called \textbf{primitive} if $A/I$ is primitive.

\medskip\textbf{Lemma 5.1.} {\it A Novikov algebra $A$ has an almost faithful irreducible module if and only if $A$ is primitive. Moreover, a primitive algebra is prime.}

\medskip\textbf{Proof.} Let $M$ be an almost faithful irreducible $A$-module. By Lemma 3.2, $M$ is isomorphic to the quotient module $A/I$ for the maximal right ideal $I$ of $A$. Let $J$ be a two-sided ideal contained in $I$. Then $MJ = 0$, which implies that $J=0$ since $M$ is almost faithful. Thus, $A$ is primitive.

\medskip Conversely, let $A$ be primitive and $I$ be the maximal modular right ideal that does not contain nonzero two-sided ideals of $A$. Then $A/I$ is an irreducible $A$-module by Lemma 3.3. Let $\rho$ be the corresponding representation, $x\in K_{\rho}(A)$. By Lemma 3.2, either $I=0$ (and $A$ is an almost faithful irreducible $A$-module), or there exist $e,f\in A$ such that $x-(x,e,f)\in I$, so that $x-(ex)f+e(xf)\in I$. Since $x\in K_{\rho}(A)$, we have $ex\in I$ and $xf\in K_{\rho}(A)$, whence $e(xf)\in I$. Thus, $x\in I$, that is, $K_{\rho}(A)\subseteq I$. Since $I$ does not contain nonzero two-sided ideals of $A$, we have $K_{\rho}(A)=0$.

By what was proved above and Lemma 3.1, the primitive algebra is prime. \qed

\medskip $\mathcal{J}$-semisimple algebras in the associative case are sometimes called \textbf{semiprimitive}. This definition is also justified for Novikov algebras as shown in the following theorem.

\medskip\textbf{Lemma 5.2.} {\it A primitive Novikov algebra is $\mathcal{J}$-semisimple.}

\medskip\textbf{Proof.} Let $A$ be a primitive Novikov algebra. By Lemma 5.1, $A$ has a maximal modular right ideal $I$ that does not contain nonzero two-sided ideals of $A$. By Lemma 3.2, two cases are possible.

\medskip 1) The module $M$ is isomorphic to $A/I$, where $I$ is a two-sided ideal, and $A/I$ is a field. Then $I=0$ and $A$ is a field.

\medskip 2) The module $M$ is isomorphic to $A/I$, where $I$ is the maximal right ideal of $A$, and $a-(a,e,f)\in I$ for any $a\in A$. The corresponding representation is almost faithful, so $I$ does not contain nonzero two-sided ideals. Suppose that $\mathcal{J}(A)\neq 0$. Then $\mathcal{J}(A)+I=A$.

\medskip Suppose $(A,A,\mathcal{J}(A))\nsubseteq I$. Then, by Lemma 3.2, we can assume that $f\in \mathcal{J}(A)$. Note that since $(a,e,f)=(e,a,f)$, then we can assume that $e\in \mathcal{J}(A)$. This means that $a-(a,e,f)\in\ I\cap \mathcal{J}(A)$ for any $a\in\mathcal
{K}(A)$. Thus, $I\cap\mathcal{J}(A)$ is a modular right ideal of $\mathcal{J}(A)$, which embeds in the maximal modular right ideal by Lemma 3.6. By Lemma 3.3, there exists an irreducible $\mathcal{J}(A)$-module, it is a contradiction.

\medskip Thus, $(A,A,\mathcal{J}(A))\subseteq I$. Then $(A,A,\mathcal{J}(A))=0$ and $\mathcal{J}(A)$ is an associative Novikov algebra. In particular, $\mathcal{J}(A)[\mathcal{J}(A),\mathcal{J}(A)]=0$ by \eqref{RAs}. By Lemma 3.1, $A$ is prime, so $\mathcal{J}(A)$ is commutative. Then $I_1=I\cap\mathcal{J}(A)$ is a two-sided ideal of $\mathcal{J}(A)$. But then $AI_1=II_1+\mathcal{J}(A)I_1 = II_1+I_1\mathcal{J}(A)\subseteq I_1+II_1$. It remains to note that $II_1\subseteq I$ and $II_1\subseteq \mathcal{J}(A)$, so $AI_1\subseteq I_1$ and $I_1$ are two-sided ideals of $A$ contained in $I$. In particular, $\mathcal{J}(A)I\subseteq \mathcal{J}(A)\cap I = 0$. Thus, $I_1=0$. But then $AI \subseteq I^2+\mathcal{J}(A)I \subseteq I$, that is, $I$ is a two-sided ideal of $A$ and $I=0$. But then $A=\mathcal{J}(A)$, it is a contradiction. \qed

\medskip To obtain the final result, we need to learn how to lift maximal modular right ideals from some ideal to the entire algebra.

\medskip\textbf{Lemma 5.3.} {\it Let $A$ be a Novikov algebra, $I$ be an ideal of $A$, and $J$ be the maximal modular right ideal of $I$. Then $J$ is a right ideal of $A$.}

\medskip\textbf{Proof.} By definition, two cases are possible.

\medskip 1) $J$ is a two-sided ideal in $I$, where $I/J$ is a field. Then $J$ is an ideal in $A$ by (\cite{P2022}, lemma~6).

\medskip 2) There exist elements $e,f\in I$ such that $a-(a,e,f)\in J$ for any $a\in I$. Suppose that $JA\nsubseteq J$. Then $Jx\nsubseteq J$ for some $x\in A$. Let $j\in J, i\in I$. Then 
\[(jx)i \overset{\eqref{RCom}}{=} (ji)x \in Jx,\] 
so $Jx$ is a right ideal of $I$. Then $J+Jx$ is a right ideal of $I$ that strictly contains $J$, that is, $J+Jx=I$. In particular, $e=c+dx$, where $c,d\in J$. Then for any $a\in I$ we have $a-(a,e,f) = a-(a,c,f)-(a,dx,f)\in J$. Note that $(a,c,f)=(c,a,f)\in J$ by \eqref{LAs}, $(a,dx,f) = (d,ax,f)\in J$ by \eqref{LAs} and \eqref{AsLin}. Thus, $a\in J$ for any $a\in I$, it is a contradiction. Hence, $JA\subseteq J$ and $J$ is a right ideal of $A$. \qed

\medskip\textbf{Lemma 5.4.} {\it Let $B$ be an ideal of a Novikov algebra $A$. Then the set of maximal modular right ideals of $B$ coincides with the set of right ideals of the form $I\cap B$, where $I$ is the maximal modular right ideal of $A$ that does not contain $B$.}

\medskip\textbf{Proof} Let $J$ be the maximal modular right ideal of $B$. By Lemma 5.3, $J$ is a right ideal of $A$. Two cases are possible.

\medskip 1) $J$ is a two-sided ideal and the algebra $B/J$ is a field. Let $e\in B$ be the preimage of the unit of the algebra $B/J$. Consider the right ideal $J_1=\{a-ea\mid a\in A\}$ and the right ideal $J_2=J_1+J$. If $e\in J_2$, then $e=a-ea+b$ for some $a\in A$, $b\in J$. Then $e^2=ea-e(ea)+eb$. Since $ea\in B$, we have $ea-e(ea)\in J$, whence $e^2=ea-e(ea)+eb\in J$, a contradiction. Thus, $e\notin J_2$. Then the set of right ideals containing $J_2$ and not containing $e$ contains a maximal element $J_0$. If $J_0\subsetneq K$, where $K$ is a right ideal, then $e\in K$. But $a - ea \in J_1 \subseteq K$ for any $a \in A$, therefore $a \in K$ and $K = A$. This means that $J_0$ is a maximal right ideal of $A$. Since $a - ea \in J_0$ for any $a \in A$ and $e \notin J_0$, the module $A/J_0$ is irreducible. From Lemma 3.2 it follows that $J_0$ is a maximal modular right ideal of $A$.

\medskip It remains to note that $B \cap J_0$ is a right ideal of $A$ satisfying the condition $J \subseteq B \cap J_0 \subsetneq B$, whence $B \cap J_0=J$ since $J$ is maximal in $B$.

\medskip 2) There exist elements $e,f\in B$ such that $a-(a,e,f)\in J$ for any $a\in B$. Consider the right ideal $J_1=\{a-(a,e,f)\mid a\in A\}$ and the right ideal $J_2=J_1+J$. If $e\in J_2$, then $e=a-(a,e,f)+b$ for some $a\in A$, $b\in J$. Then $e^2 = ae-(ae,e,f)+be$. Since $ae\in B$, we have $ae-(ae,e,f)\in J$. Since $J$ is a right ideal, we have $be\in J$. Thus, $e^2\in J$, it is a contradiction. Thus, $e\notin J_2$. Then the set of right ideals containing $J_2$ and not containing $e$ contains a maximal element $J_0$. As in the previous case, $J_0$ is a maximal right ideal of $A$. Since $a-(a,e,f)\in J_0$ for any $a\in A$, $J_0$ is a maximal modular right ideal of $A$. It remains to note that $B\cap J_0$ is a right ideal of $A$ satisfying the condition $J\subseteq B\cap J_0\subsetneq B$, whence $B\cap J_0=J$ since $J$ is maximal in $B$. \qed

\medskip In some classical cases (associative, alternative), the Jacobson radical coincides with the intersection of the kernels of all irreducible representations. In the case of Novikov algebras, a different, slightly weaker result holds.

\medskip\textbf{Theorem 5.5.} {\it The Jacobson radical of a Novikov algebra coincides with the intersection of the quasi-kernels of all its irreducible representations.}

\medskip\textbf{Proof.} Denote by $L$ the intersection of $K_{\rho}(A)$ over all irreducible representations $\rho$ of the Novikov algebra $A$. For any such $\rho$, by Lemma 5.1, $\mathcal{J}(A/K_{\rho}(A))=0$. Therefore, $\mathcal{J}(A)\subseteq K_{\rho}(A)$ by Theorem 4.2 due to the well-known properties of the radical. Thus, $\mathcal{J}(A)\subseteq L$.

Suppose that $\mathcal{J}(A)\subsetneq L$. Then there exists an irreducible $L$-module. By Lemma 3.2, there exists a maximal right modular ideal $I$ of $L$. By Lemma 5.4, $I=J\cap L$ for some maximal modular right ideal $J$ of $A$, it is a contradiction. Thus, $\mathcal{J}(A)=L$. \qed

\medskip From Example 3.4 it follows that the intersection of kernels of irreducible representations need not even be a subalgebra, much less coincide with $\mathcal{J}(A)$, even in the finite-dimensional simple algebra.

\medskip By Lemma 5.1, there is a one-to-one correspondence between the primitive ideals of $A$ and the ideals of $K_{\rho}(A)$ with respect to the irreducible representations of $A$. Thus, the following corollary holds.

\medskip\textbf{Corollary 5.6.} {\it In every Novikov algebra, the Jacobson radical coincides with the intersection of all its primitive ideals. Every  $\mathcal{J}$-semisimple Novikov algebra is a subdirect product of primitive algebras.}

\medskip Note that in the finite-dimensional case, the Jacobson radical is a solvable ideal, as the following assertion shows.

\medskip\textbf{Proposition 5.7.} {\it Let $A$ be a finite-dimensional Novikov algebra. Then $\mathcal{J}(A)$ coincides with the largest solvable ideal $N(A)$.}

\medskip\textbf{Proof.} Let $C$ be a simple finite-dimensional Novikov algebra and $I$ be the maximal right ideal of $C$. If $(C/I)C=0+I$, then $C^2\subseteq I$, that is, $C^2\neq C$. But then $C^2=0$ since $C$ is simple, therefore $(C/I)C\neq 0$, which means the module $C/I$ is irreducible. Thus, a simple finite-dimensional Novikov algebra has an irreducible module, and hence a classically semisimple nonzero finite-dimensional Novikov algebra also has an irreducible module.

Let $A$ be a finite-dimensional Novikov algebra. From Lemma 3.2, it is clear that $N(A)$ has no irreducible right modules, hence $N(A)\subseteq\mathcal{J}(A)$. Consider the algebra $B=\mathcal{J}(A)/N(A)$. The algebra $B$ is a semiprime $\mathcal{J}$-radical Novikov algebra. Then $B$ is classically semisimple, and hence (if $B\neq 0$) it has an irreducible right module. Since $B$ is $\mathcal{J}$-radical, this means that $B=0$. \qed

\medskip Novikov algebras are not monoassociative, i.e. one-generated Novikov does not have to be associative. It means that there are different notions of element nilpotency. We will need one of them which earlier appeared in \cite{DT2006, P2026}.

\medskip\textbf{Definition.} An element $x\in A$ is called \textbf{r-nilpotent} if $x^n = 0$ for some $n>0$. A Novikov algebra $A$ is called {\bf r-nil} if every element of $A$ is r-nilpotent. An ideal $I$ of $A$ is called {\bf r-nil} if every element of $I$ is r-nilpotent.

\medskip Corollary~5.6 implies that Jacobson radical of Novikov algebra contains Baer radical. One of the important properties of the Jacobson radical of associative algebras is that it contains all nil ideals. Let us prove the similar statement for Novikov algebras.

\medskip\textbf{Proposition 5.8.} {\it Let $A$ be a Novikov algebra with an r-nil nonzero ideal $L$. Then $L\subseteq J(A)$.}

\medskip\textbf{Proof.} Suppose that $L\nsubseteq J(A)$. Then an image of $L$ under canonical projection $A\to A/J(A)$ is a nonzero $J$-semisimple r-nil ideal. By Corollary~5.6 there exists a nonzero primitive r-nil algebra $B$. By Lemma~5.1 there exists an almost faithful irreducible $B$-module $M$. By Lemma~3.2 there exists a maximal modular right ideal $I$ of $B$ such that $M\simeq B/I$. Since the algebra $B$ is r-nil, then $B/I$ is not a field. It means that there exist elements $e,f\in B$ such that $x-(x,e,f)\in I$ for any $x\in B$. The element $e\in A$ is r-nilpotent, $e^n=0$. Hence, 
\[(x-(x,e,f))R_e^{n-1} \overset{\eqref{AsLin}}= xR_e^{n-1}\in I.\]
Then $MR_e^{n-1} = 0$ and $Me=0$. Therefore $x+x(ef)\in I$ for any $x\in B$. 
Denote $g:=ef$. Then $x+xg\in I$ for any $x\in B$. Hence, $xR_g^k + xR_g^{k+1}\in I$ for any $x\in B$. Therefore, $x+(-1)^{k+1}xR_g^{k}\in I$ for any integer $k\ge 1$. 

\medskip Let $g^m=0$. We will prove the following statement:
\[xR_g^m = (xg^k)R_g^{m-k}\]
for any integer $k\ge 1$. For $k=1$ the statement is obvious. If $k>1$ then 
\begin{gather*}
	xR_g^m = (xg^k)R_g^{m-k} = (x,g^k,g)R_g^{m-k-1} + xg^{k+1}R_g^{m-k-1} \overset{\eqref{AsLin}}{=} \\ \overset{\eqref{AsLin}}{=} (x,g^{m-1},g)+xg^{k+1}R_g^{m-k-1} = xg^{k+1}R_g^{m-k-1}.
\end{gather*}
If $k=m$ then we have $xR_g^m = xg^m = 0$. Therefore, $x+(-1)^{m+1}xR_g^m = x$ and $x\in I$ for any $x\in B$, a contradiction. \qed

\medskip\textbf{Definition.} A radical class $\mathcal{R}$ is called \textbf{hereditary} if for any algebra $A$ and any of its ideals~$I$, the equality
\[\mathcal{R}(I)=\mathcal{R}(A)\cap I\]
holds.

\medskip All classical radicals in associative, alternative, or Jordan algebras are usually hereditary. The Jacobson radical of Novikov algebras is no exception.

\medskip\textbf{Theorem 5.9.} {\it The Jacobson radical of Novikov algebras is hereditary.}

\medskip\textbf{Proof.} It suffices to prove that the ideal of a $\mathcal{J}$-semisimple algebra is a $\mathcal{J}$-semisimple algebra, and the ideal of a $\mathcal{J}$-radical algebra is a $\mathcal{J}$-radical algebra.

\medskip Let $A$ be a $\mathcal{J}$-semisimple algebra and $I$ be an ideal of $A$. Since $I\neq\mathcal{J}(I)$, the algebra $I/\mathcal{J}(I)$ is $\mathcal{J}$-semisimple. It is easy to see that an algebra with zero multiplication is $\mathcal{J}$-radical, so the algebra $I/\mathcal{J}(I)$ is semiprime. Then $\mathcal{J}(I)$ is a $\mathcal{J}$-ideal of $A$ (by \cite{P2022}), that is, $\mathcal{J}(I)=0$.

\medskip Let $A$ be a Novikov algebra and $I$ be an ideal of $A$. Suppose that $I\neq\mathcal{J}(I)$. Then the algebra $I$ has an irreducible module $M$. By Lemma 3.2, $M$ is isomorphic to the quotient module $I/J$ for some right maximal modular ideal $J$ of $I$. By Lemma 5.4, the algebra $A$ contains a maximal modular right ideal. By Lemma 3.3, this means that there exists an irreducible $A$-module, that is, $A\neq\mathcal{J}(A)$. \qed

\section{Left Novikov Representations}

All definitions from the previous sections can be adapted for left representations and left modules. In this small section, we will see through examples that left representations of Novikov algebras do not share all the desirable properties established for right representations earlier. A single example due to Zel'manov suffices.

\medskip\textbf{Example 6.1.} (\cite{Z1987}) Let $A$ be a two-dimensional algebra, $A=L(a,b)$ with the multiplication table $a^2=b^2=ba=0$, $ab=b$. Then $A$ is a solvable Novikov algebra with $(A^2)^2=0$. The element $a$ generates a left ideal $I=L(a)$. Consequently, $A/I$ is an irreducible left $A$-module since $a(b+I) = b+I\neq 0+I$. Therefore, a solvable Novikov algebra can admit an irreducible left representation, which is not the case for right representations, by Corollary 5.6.

\medskip Let $J=L(b)$. Then $J$ is a two-sided ideal of $A$. Since $J^2=0$, $J$ does not admit any irreducible left representation. Furthermore, since $(A/J)^2=0$, $A/J$ does not admit irreducible left representation either. Thus, the class of Novikov algebras that do not admit irreducible left representations is not closed under extensions. Note, however, that this property fails for right representations by Lemma 4.1. Therefore, an analogue of the left Jacobson radical does not exist for Novikov algebras.

\medskip Evidently, these observations are enough to refrain from a deeper investigation of left representations, given their inherent limitations.

\section*{Acknowledgments}

The author was supported by RAS Fundamental Research Program (project FWNF-2026-0017).

\noindent Alexander Panasenko \\
Sobolev Institute of Mathematics \\
Acad. Koptyug ave. 4, 630090 Novosibirsk, Russia \\
e-mail: a.panasenko@g.nsu.ru


\begin{thebibliography}{99}

\bibitem{GD1979} I.M. Gel'fand, I.Ya. Dorfman, Hamiltonian operators and algebraic structures related to them, Funct. Anal. Appl. 13 (4) (1979) 248–262.

\bibitem{BN1985} A.A. Balinskii, S.P. Novikov, Poisson brackets of hydrodynamic type, Frobenius algebras and Lie algebras, Sov. Math. Dokl. 283 (5) (1985) 1036--1039.

\bibitem{N1985} S.P. Novikov, The geometry of conservative systems of hydrodynamic type. The method of averaging for field-theoretical systems, Russian Math. Surveys  40 (4) (1985) 85--98.

\bibitem{Z1987} E.I. Zel'manov, A class of local translation-invariant Lie algebras, Sov. Math. Dokl 292 (6) (1987) 1294--1297.

\bibitem{F1989} V.T. Filippov, A class of simple nonassociative algebras, Math. Notes 45 (1) (1989) 68--71.

\bibitem{Osborn1992} J.M. Osborn, Simple novikov algebras with an idempotent, Communications in Algebra 20 (9) (1992) 2729--2753.

\bibitem{Osborn1995} J.M. Osborn, Modules over Novikov algebras of characteristic 0, Communications in Algebra 23 (10) (1995) 3627--3640.

\bibitem{OZ1995} J.M. Osborn, E. Zelmanov, Nonassociative algebras related to Hamiltonian operators in
the formal calculus of variations, Journal of Pure and Applied Algebra 101 (1995) 335--352.

\bibitem{Xu1996} X. Xu, On Simple Novikov Algebras and Their Irreducible Modules, Journal of Algebra 185 (1996) 905--934.

\bibitem{ZZ2024} V.N. Zhelyabin, A.S. Zakharov, On finite-dimensional simple Novikov algebras of characteristic p, Siber. Math. J. 65 (3) (2024) 680--687.

\bibitem{Xu2001} X. Xu, Classification of Simple Novikov Algebras and Their Irreducible Modules of Characteristic 0, Journal of Algebra 246 (2001), 673--707.

\bibitem{SS1974} A.M. Slin'ko, I.P. Shestakov, Right representations of algebras, Algebra Logika 13(5) (1974), 544--588.

\bibitem{P2022} A.S. Panasenko, Semiprime Novikov algebras, Int. J. of Algebra and Comp. 32 (7) (2022) 1369--1378.

\bibitem{P2024} A.S. Panasenko, On radicals of Novikov algebras, Comm. in Algebra 52 (1) (2024) 140--147.

\bibitem{PZ2026} A. Pozhidaev, V. Zhelyabin, On simple and semisimple finite-dimensional Novikov algebras and their automorphisms, Journal of Algebra 689 (2026) 1--26.

\bibitem{SF2002} V.A. Sereda, V.T. Filippov, On Homotopes of Novikov Algebras, Siberian Math. J. 43 (1) (2002) 1--7.

\bibitem{P2026} A.S. Panasenko, Radicals of Lie-solvable Novikov algebras, Int. J. of Algebra and Comp., online ready DOI: 10.1142/S0218196726500451.

\bibitem{SZ2020} I. Shestakov, Z. Zhang, Solvability and nilpotency of Novikov algebras, Comm. in Algebra 48 (2020) 5412--5420.

\bibitem{DT2006} A.S. Dzhumadil'daev, K.M. Tulenbaev, Engel Theorem for Novikov Algebras, Comm. in Algebra 34 (2006) 883--888.

\end{thebibliography}
\end{document}